\documentclass[%
 reprint,
 amsmath,amssymb,
 aps,
]{revtex4-2}
\usepackage{graphicx}
\newcommand{\physrep}{Phy. Report}
\usepackage{dcolumn}
\usepackage{bm}
\usepackage{tikz}
\usetikzlibrary{angles}
\usetikzlibrary{quotes}
\usepackage[utf8]{inputenc}
\usepackage[T1]{fontenc}
\usepackage{mathptmx}
\usepackage{amsmath,amssymb}
\usepackage{pifont}
\usepackage{graphicx}
\usepackage{dcolumn}
\usepackage{bm}

\begin{document}

\preprint{APS/123-QED}
\title{Orbits of charged particles with an azimuthal initial velocity in a dipole magnetic field } 



\author{Hanrui Pang}
\affiliation{ School of Physical Science and Technology, Southwest Jiaotong University, Chengdu 610031, People’s Republic of China}
\author{Siming Liu}\email{liusm@swjtu.edu.cn}
\affiliation{ School of Physical Science and Technology, Southwest Jiaotong University, Chengdu 610031, People’s Republic of China}
\author{Rong Liu}
\affiliation{ School of Physical Science and Technology, Southwest Jiaotong University, Chengdu 610031, People’s Republic of China}




%
\date{\today}

\begin{abstract}
Nonintegrable dynamical systems have complex structures in their phase space. 
Motion of a test charged particle in a dipole magnetic field can be reduced to a 2 degree-of-freedom (2 d.o.f.) nonintegrable Hamiltonian system. We carried out a systematic study of orbits of charged particles with an azimuthal initial velocity in a dipole field via 
calculation of their Lyapunov characteristic exponents (LCEs) and escape times for a dimensionless energy less and greater than 1/32, respectively. Meridian plane periodic orbits symmetric with respect to the equatorial plane are then identified. We found that 1) symmetric periodic orbits can be classified into several classes based on their number of crossing points on the equatorial plane; 2) the initial conditions of these classes locate on closed loops or closed curves going through the origin; 3) most isolated regions of stable quasi-periodic orbits are associated asymmetric stable periodic orbits; 4) classes of asymmetric periodic orbits either go through the origin or terminate at flat equatorial plane orbits with the other end approaching centers of spiral structures; 5) there are apparent self-similarities in the above features with the decrease of energy.

\end{abstract}

\maketitle 

\section{\label{sec:level1} Introduction 
}

The behavior of nonintegrable dynamical systems is complicated due to presence of complex structures in their phase space \citep{1979PhR....52..263C}.
In particular, stable trajectories wherever exist can affect properties of the whole system significantly \citep{PhysRevE.97.022215}. Although chaotic trajectories are ergodic, they may also show nontrivial behaviors \citep{lichtenberg_1992, 2002PhR...371..461Z}. 
As one of the simplest nonintegrable Hamiltonian systems, motion of charged particles in the dipole magnetic field has been studied extensively during the past century \citep{1930ZA......1..237S, 1965articleD, 2022Chaos..32d3104L}. In particular, significant efforts have been put to study periodic orbits in the Meridian plane \citep{devogelaere1950, 1977CeMec..51...177, 1977Ap&SS..48..471M, 1978CeMec..17..215M, 1978CeMec..17..233M}, and it has been known that stable periodic orbits can affect statistical properties of the whole system significantly \citep{2002PhR...371..461Z, PhysRevE.97.022215}. Stable periodic orbits also bear significance in physical systems since their associated stable quasi-periodic orbits can be realized and have distinct characteristics \citep{2021ApJ...909...59W, 2018PhRvL.121w5003H}.

With the aid of Lyapunov characteristic exponents (LCE), stable orbits can be readily identified since the value of their maximum LCE (mLCE) is zero. However, the distribution of stable orbits in the phase space is complex, showing fractal structures \citep{1985PhRvL..55..661U, SeokLee_2020, 2022Chaos..32d3104L}. To understand this distribution, one may identify classes of periodic orbits \citep{1978CeMec..17..233M, 1987PhyD...26..369P}. In this paper, we focus on studying orbits with an azimuthal initial velocity, i.e., initially static in the Meridian plane. We first show that periodic orbits in the Meridian plane have two static points half-a-period apart and those symmetric with respect to the equatorial plane cross it perpendicularly twice in each period at the same distance from the origin. Classes of such symmetric open periodic orbits are then identified and used to explore distribution of stable orbits in the phase space.

In \S \ref{stormer} we present the Hamiltonian equations for the Meridian plane motion of a charged particle moving in a dipole magnetic field and prove general properties of periodic orbits with open ends in \S \ref{proof}. Our main results are shown in \S \ref{results}. Conclusions are drawn in \S \ref{conc}.

\section{Motion of charged particles in the Meridian plane of a dipole magnetic field}  
\label{stormer}

Since the dipole magnetic field is axisymmetric, the azimuthal canonical momentum $p_\phi$ and the energy $H$ of a charged particle moving in it are conserved if one ignores the radiation reaction force, which is a quantum phenomenon and cannot be treated self-consistently in classical mechanism \citep{2022PhRvD.105a6024P, 2009PhRvD..80b4031G}. Following the standard procedure \citep{1965articleD, 2022Chaos..32d3104L}, one can introduce a characteristic length with $p_\phi$, the charge of the particle $e$, and magnetic moment of the dipole $M$, and obtain the dimensionless Hamiltonian for motion in the Meridian plane $[\rho, z]$:
\begin{equation}
    H= \frac{1}{2}\left({p_z}^2+{p_\rho}^2\right)+ V ,
\end{equation}
where $z$ is the symmetric axis, $\rho$ is the radius in the equatorial plane, and the effective potential
\begin{equation}
    V=\frac{1}{2}\left[\frac{1}{\rho}-\frac{\rho}{(z^2+\rho^2)\,^\frac{3}{2}}\right]^2
\end{equation}   
corresponds to the kinetic energy associated with the azimuthal motion, $p_\rho$, and $p_z$ are the momentum in the $\rho$, and $z$ directions, respectively. 
We then have the Hamiltonian differential equations of motion:
       \begin{align}
         \dot{z}&= \frac{\partial H}{\partial p_z}=p_z\,, \\[1ex]
  \dot{\rho}&=\frac{\partial H}{\partial p_\rho}=p_\rho\,, \\[1ex]
  \dot{p_z}&= -\frac{\partial H}{\partial z}=-\frac{3z\rho\left[\frac{1}{\rho}-\frac{\rho}{(z^2+\rho^2)^\frac{3}{2}}\right]}{(z^2+\rho^2)^\frac{5}{2}}\,,  \\
      \dot{p_\rho} &=-\frac{\partial H}{\partial \rho}\ ,\\[1ex] \nonumber
  &=\left[ \frac{1}{\rho^2}-\frac{3\rho^2}{(z^2+\rho^2)^\frac{5}{2}}
  +\frac{1}{(z^2+\rho^2)^\frac{3}{2}}\right]\left[\frac{1}{\rho}-\frac{\rho}{(z^2+\rho^2)^\frac{3}{2}}\right],
  \end{align}
where an upper dot $\dot{}$ indicates a derivative with respect to the time $t$.

\section{General Characteristics of Open Periodic Orbits}
\label{proof}

Periodic orbits are the simplest orbits and are ubiquitous in the phase space \citep{https://doi.org/10.1002/cpa.3160260204}. Stable quasi-periodic orbits are always associated with stable periodic orbits \citep{1987PhyD...26..369P, lichtenberg_1992}. It is therefore essential to understand properties of periodic orbits to uncover the phase space structure of a dynamical system. 

Particles with an azimuthal initial velocity have a zero initial velocity with $p_{z,0}=0,p_{\rho,0}=0$ in the $[z,\rho]$ plane. Periodic orbits in this Meridian plane have open paths and some general characteristics that are essential to search for such orbits numerically. To prove these characteristics, we define $\Ddot{z}\equiv f(z,\rho)$ and $\Ddot{\rho}\equiv g(z,\rho)$, and the minimum positive period $T$. Static orbits with $T=0$ are trivial and will not be considered here.
The Hamiltonian equations then lead to
\begin{equation}
\label{e-1}
\begin{split}
    f(-z,\rho) &= -f(z,\rho)\,,\\
    g(-z,\rho) &= g(z,\rho)\,,\\
\end{split}
\end{equation}
and
\begin{equation}
\label{e0}
    \begin{split}
    z(t+T)=z(t),\ \ \dot{z}(t+T)=\dot{z}(t)\,,\\
    \rho(t+T)=\rho(t),\ \ \dot{\rho}(t+T)=\dot{\rho}(t)\,.\\
    \end{split}
\end{equation}


Since open orbits start with zero initial velocities, it is evident that $z(t)=z(-t),\ \rho(t)=\rho(-t)$. Let $t=-t^*$, we have 
\begin{equation}
\label{e2}
    \begin{cases}
    \frac{d^2z}{{dt^*}^2}=f(z,\rho)\\
    \frac{d^2\rho}{{dt^*}^2}=g(z,\rho)  
    \end{cases}
    {\rm with}\ \ 
    \begin{cases}
    z|_{t^*=0}=z_0,\ \ \frac{dz}{dt^*}|_{t^*=0}=0\,,\\
    \rho|_{t^*=0}=\rho_0,\ \ \frac{d\rho}{dt^*}|_{t^*=0}=0 \,,
    \end{cases}
\end{equation}
which are the same equations followed by $z(t)$ and $\rho(t)$.
Then we have $z(t)=z(t^*),\rho(t)=\rho(t^*)$, i.e., 
\begin{equation}
\label{e3}
\begin{split}
    z(t)=z(-t)\,,\\
    \rho(t)=\rho(-t)\,.
\end{split}
\end{equation} 

\subsection{Open Periodic Orbits}

For open periodic orbits, it can be shown that there are two static points in each period at times $t=0$ and $t=\frac{T}{2}$, respectively.
From equation (\ref{e3}), we have
\begin{equation}
\label{e3.1}
\begin{split}
    \dot{z}(t)=-\dot{z}(-t)\,,\\
    \dot{\rho}(t)=-\dot{\rho}(-t)\,.
\end{split}
\end{equation} 
Let $t=\frac{T}{2}$, equation (\ref{e0}) then gives 
\begin{equation}
    \begin{split}
    \dot{z}(\frac{T}{2})=-\dot{z}(-\frac{T}{2})=-\dot{z}(\frac{T}{2})\,,\\
    \dot{\rho}(\frac{T}{2})=\dot{\rho}(-\frac{T}{2})=-\dot{\rho}(\frac{T}{2})\,.
    \end{split}
\end{equation}
Therefore
\begin{equation}
\label{e5}
    \begin{split}
    \dot{z}(\frac{T}{2})=0\,,\\
    \dot{\rho}(\frac{T}{2})=0\,,
    \end{split}
\end{equation}
and we define $z_1\equiv z(\frac{T}{2}),\ \rho_1\equiv \rho(\frac{T}{2})$. Clearly, $(z_1,\rho_1)\neq(z_0,\rho_0)$. Otherwise, $T$ will no longer be the minimum positive period since $T/2< T$ is also a positive period.


Furthermore, we can prove that $t=\frac{T}{2}$ is the only $t$ that satisfies $\dot{z}(t)=0,\ \dot{\rho}(t)=0$ in $(0,T)$.
Suppose that there exists a $t_1\in(0,T),t_1\neq\frac{T}{2}$ such that $\dot{z}(t_1)=0,\ \dot{\rho}(t_1)=0$,
we have
\begin{equation}
    \label{e7}
    \begin{cases}
    \frac{d^2z}{{dt}^2}=f(z,\rho)\\
    \frac{d^2\rho}{{dt}^2}=g(z,\rho)  
    \end{cases}
    {\rm with}\ \ 
    \begin{cases}
    z|_{t=t_1}=z(t_1),\ \ \frac{dz}{dt}|_{t=t_1}=0\,,\\
    \rho|_{t=t_1}=\rho(t_1),\ \ \frac{d\rho}{dt}|_{t=t_1}=0\,. 
    \end{cases}
\end{equation}
Let $t=-t^{*}+2t_1$, we have
\begin{equation}
    \label{e7.1}
    \begin{cases}
    \frac{d^2z}{{dt^*}^2}=f(z,\rho)\\
    \frac{d^2\rho}{{dt^*}^2}=g(z,\rho)  
    \end{cases}
    {\rm with}\ \ 
    \begin{cases}
    z|_{t^*=t_1}=z(t_1),\ \ \frac{dz}{dt^*}|_{t^*=t_1}=0\,,\\
    \rho|_{t^*=t_1}=\rho(t_1),\ \ \frac{d\rho}{dt^*}|_{t^*=t_1}=0 \,.
    \end{cases}
\end{equation}
Equations (\ref{e7}) and (\ref{e7.1}) have the same solutions: $z(t)=z(t^*),\ \rho(t)=\rho(t^*)$. Therefore
\begin{equation}
    \begin{split}
    z(t)=z(-t+2t_1)=z(t-2t_1)=z(t+T-2t_1)\,,
    \\
    \rho(t)=\rho(-t+2t_1)=\rho(t-2t_1)=\rho(t+T-2t_1)\,.
    \end{split}
\end{equation}
Then $0<|T-2t_1| \ {\rm or}\ 2t_1<T$ is a positive period. 

\subsection{Symmetric Open Periodic Orbits}
\label{Symmetric Orbits}
For open periodic orbits symmetric with respect to the equatorial plane, $(z_0,\rho_0) = (-z_1,\rho_1)$. When $t=\frac{T}{4}$, it can be shown that they pass through the equatorial plane perpendicularly with $p_\rho=0$ and $z=0$ and there is only one such perpendicular crossing within $\frac{T}{2}$.  Let $z=-z^*,\ t=t^*+\frac{T}{2}$ and considering equation (\ref{e-1}), we have 
\begin{equation}
    \label{e4}
    \begin{cases}
    \frac{d^2z^*}{{dt^*}^2}=f(z^*,\rho)\\
    \frac{d^2\rho}{{dt^*}^2}=g(z^*,\rho)  
    \end{cases}
    {\rm with} \ \
    \begin{cases}
    z^*|_{t^*=0}=-z_1=z_0,\ \ \frac{dz^*}{dt^*}|_{t^*=0}=0\,,\\
    \rho|_{t^*=0}=\rho_1=\rho_0,\ \ \frac{d\rho}{dt^*}|_{t^*=0}=0\,,
    \end{cases}
\end{equation}
which are the same equations followed by $(z(t),\ \rho(t))$.
Therefore $z(t)=z^*(t^*),\ \rho(t)=\rho(t^*)$, and we have
\begin{equation}
\label{e3.2}
    \begin{split}
    z(t)=-z(t-\frac{T}{2})\,,\\
    \rho(t)=\rho(t-\frac{T}{2})\,.
    \end{split}
\end{equation}
Taking the time derivative of $\rho(t)$, one has
\begin{equation}
\label{e3.3}
    \begin{split}
    \dot{\rho}(t)=\dot{\rho}(t-\frac{T}{2})\,.
    \end{split}
\end{equation}
Let $t=\frac{T}{4}$,  equations (\ref{e3}), (\ref{e3.1}), (\ref{e3.2}), and (\ref{e3.3}) give
\begin{equation}
    \begin{split}
    z(\frac{T}{4})=z(-\frac{T}{4})=-z(-\frac{3T}{4})=-z(\frac{T}{4})\,,\\
    \dot{\rho}(\frac{T}{4})=\dot{\rho}(-\frac{T}{4})=-\dot{\rho}(\frac{T}{4})\,.
    \end{split}
\end{equation}
Therefore 
\begin{equation}
\label{e6}
    \begin{split}
    z(\frac{T}{4})=0\,,\\
    \dot{\rho}(\frac{T}{4})=0\,.
    \end{split}
\end{equation}
One can prove that $t=\frac{T}{4}$ is the only $t\in [0,\frac{T}{2})$ that satisfies $z(t)=0,\ \dot{\rho}(t)=0$.
Suppose there exists a $t_1\in[0,\frac{T}{2})$, and $t_1\neq\frac{T}{4}$ such that $z(t_1)=0,\ \dot{\rho}(t_1)=0$, we have
\begin{equation}
    \label{e8}
    \begin{cases}
    \frac{d^2z}{{dt}^2}=f(z,\rho)\\
    \frac{d^2\rho}{{dt}^2}=g(z,\rho)  
    \end{cases}
    {\rm with}\ \ 
    \begin{cases}
    z|_{t=t_1}=0,\ \ \frac{dz}{dt}|_{t=t_1}=v_{z}(t_1)\,,\\
    \rho|_{t=t_1}=\rho_{t_1},\ \ \frac{d\rho}{dt}|_{t=t_1}=0\,. 
    \end{cases}
\end{equation}
Let $t=-t^{*}+2t_1,\ z=-z^*$, equation (\ref{e-1}) implies:
\begin{equation}
    \label{e8.1}
    \begin{cases}
    \frac{d^2z^*}{{dt^*}^2}=f(z^*,\rho)\\
    \frac{d^2\rho}{{dt^*}^2}=g(z^*,\rho)  
    \end{cases}
    {\rm with}\ \ 
    \begin{cases}
    z^*|_{t^*=t_1}=0,\ \ \frac{dz^*}{dt^*}|_{t^*=t_1}=v_{z}(t_1)\,,\\
    \rho|_{t^*=t_1}=\rho_{t_1},\ \ \frac{d\rho}{dt^*}|_{t^*=t_1}=0 \,.
    \end{cases}
\end{equation}
Equations (\ref{e8}) and (\ref{e8.1}) have the same solutions: $z(t)=z^*(t^*),\ \rho(t)=\rho(t^*)$. Equations (\ref{e3}), (\ref{e3.2}), and (\ref{e3.3}) then lead to
\begin{equation}
    \begin{split}
    z(t)=-z(-t+2t_1)=-z(t-2t_1)=z(t+\frac{T}{2}-2t_1)\,,\\
    \rho(t)=\rho(-t+2t_1)=\rho(t-2t_1)=\rho(t+\frac{T}{2}-2t_1)\,.
    \end{split}
\end{equation}
Then $0<|\frac{T}{2}-2t_1|\leq\frac{T}{2}<T$ is a  positive period. 

One can actually show that any open orbit that crosses the equatorial plane perpendicularly must be symmetric open periodic orbits.
Equations (\ref{e3}), (\ref{e8}), and (\ref{e8.1}) imply
\begin{equation}
    \begin{split}
    z(t)=-z(-t+2t_1)=-z(t-2t_1)\,,\\
    \rho(t)=\rho(-t+2t_1)=\rho(t-2t_1)\,.
    \end{split}
\end{equation}
Then we have 
\begin{equation}
    \begin{split}
    z(t)=z(t-4t_1)\,,\\
    \rho(t)=\rho(t-4t_1)\,.
    \end{split}
\end{equation}
Therefore $T=4t_1$ is a positive period and $z(0)=-z(T/2)$, $\rho(0)=\rho(T/2)$.

\section{Results}
\label{results}

\begin{figure*}[htbp]
\includegraphics[width=1.6\columnwidth, angle=0.0]{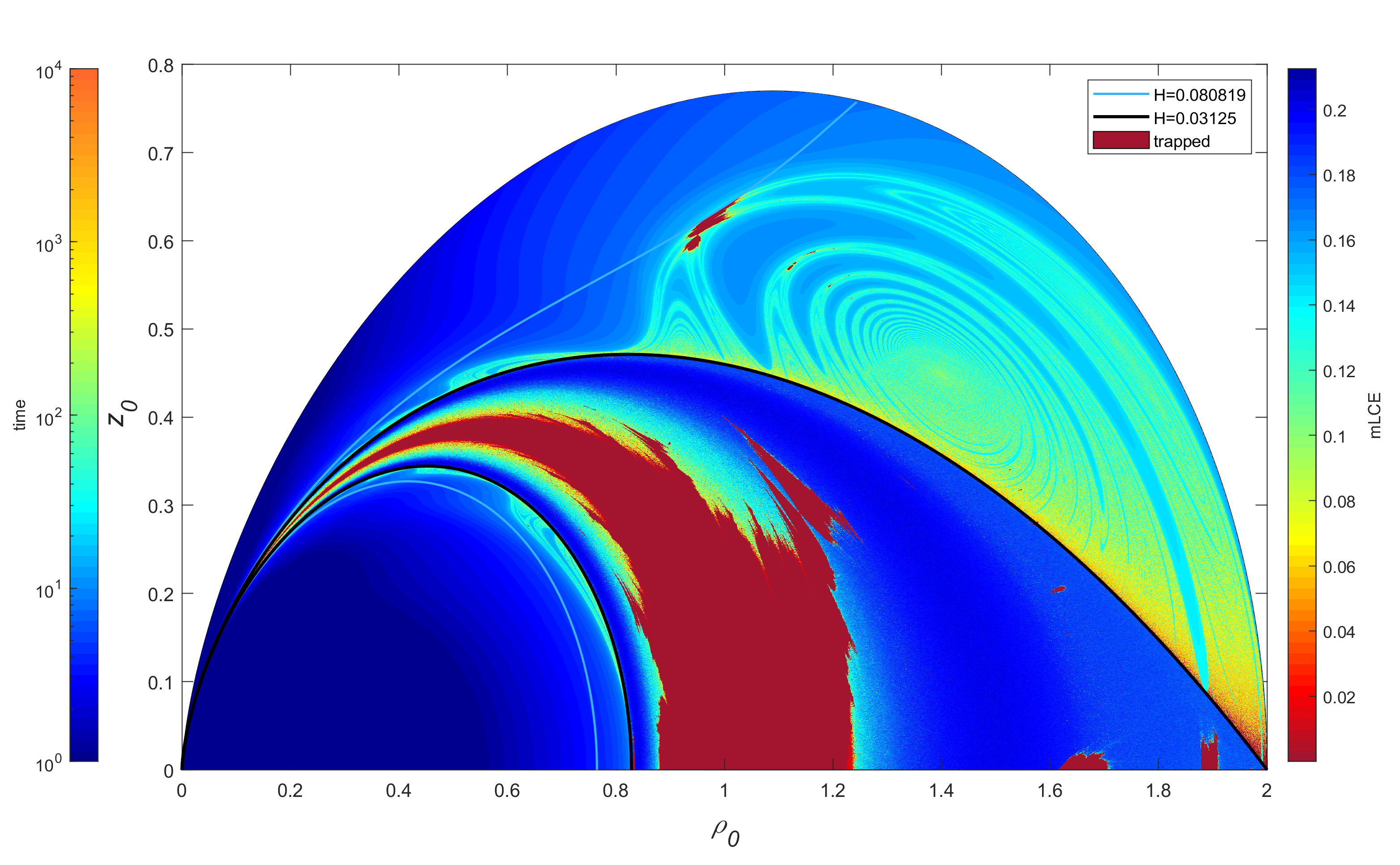}
\caption{\label{fig1} The distribution of the mLCE and the escape time of particles with an azimuthal initial velocity. The cyan lines indicate the maximum energy of periodic orbits \citep{1978CeMec..17..215M}. The radial acceleration is positive beyond the outer boundary. Besides several blocks of stable quasi-periodic orbits associated with stable periodic orbits in the equatorial plane, there are many small isolated regions of stable quasi-periodic orbits. See text for details.}
\end{figure*}


Using the same numerical codes used by \citet{2022Chaos..32d3104L}, i.e. RK45 formula being used to solve the relevant ordinary differential equations with default absolute and relative error tolerance of $10^{-9}$, we have the main result shown in Figure \ref{fig1}. For $H<1/32$, the $(\rho,\ z)$ plane is scanned to obtain the mLCE of each orbit starting with $\rho=\rho_0,\ z=z_0,\ p_\rho=p_z=0$. It should be pointed out that at a given energy, most chaotic orbits should have the same value of the mLCE \citep{SeokLee_2020}. Finite integration time and stickiness of orbits near stable orbits make the numerical results not as uniform as expected \citep{meiss_2017,2018PhRvE..97b2215H, 2019PhRvE..99c2203H}.

Since $\dot{r}>0$ with $r = (\rho^2+z^2)^{1/2}$ for $r>(4H)^{-1/3}$\,,\citep{1978CeMec..17..215M} we only consider the domain with $r\le (4H)^{-1/3}$, beyond which all particles escape to infinity. In the Meridian plane, this upper limit is equivalent to the magnetic field line passing through the point ($z=0, \ \rho=2$), then we have $r\le 2\cos^2\lambda$, where $\cos\lambda = \rho/r$. For $H\ge 1/32$, we calculate the escape time $t_{\rm esc}$ that it takes for the particle to reach $\rho=2$. Particles with $t_{\rm esc}>10^4$ are considered as trapped particles, whose orbits are likely stable. It is interesting to note that there are spiral structures in the escape time and most isolated regions of stable orbits locate at some turning points of spiral structures. Such structures also appear to emerge repeatedly approaching the origin. We also indicate the contours for the maximum energy $H\simeq0.80819$ that periodic orbits can have \citep{1978CeMec..17..215M} and $H=1/32$ in Figure \ref{fig1}.

The distribution of stable orbits (dark red) shows complex structures \citep{Xie2020FromPT}. Besides stable orbits at low energies and those associated stable orbits in the equatorial plane \citep{devogelaere1950, 1978CeMec..17..233M}, there are many small islands of stable orbits spread over the whole phase space. The boundaries of regions of stable orbits appear to be fractal, and there are apparent self-similarities with the decrease of energy.
We will try to understand these structures by identifying classes of symmetric open periodic orbits and analyzing their relations to regions of stable orbits. 

 \begin{figure*}[htb]
\includegraphics[width=1.6\columnwidth, angle=0.0]{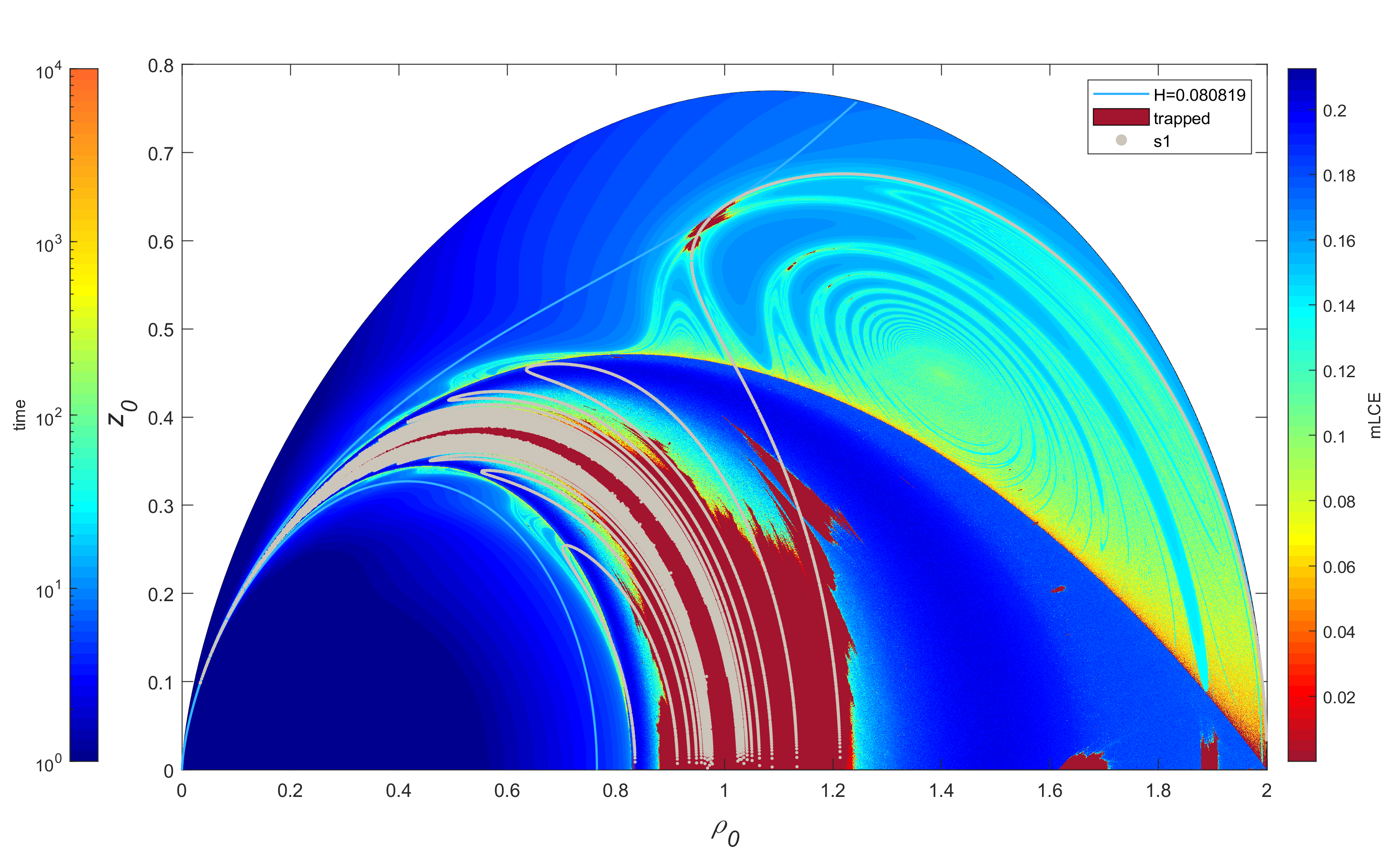}
\caption{\label{fig2} Same as Fig. \ref{fig1} with the open-path periodic orbits s1 indicated by gray lines. From left and right to middle, we have the families $f_2$, $f_6$,... and $f_0$, $f_4$, ... defined by \citet{1978CeMec..17..215M}, respectively. Regions of stable orbits at the low energy end of each family have similar structures.}
\end{figure*}

\subsection{Classes of Symmetric Open Periodic Orbits}

\begin{figure*}[htbp]
\includegraphics[width=1.6\columnwidth, angle=0.0]{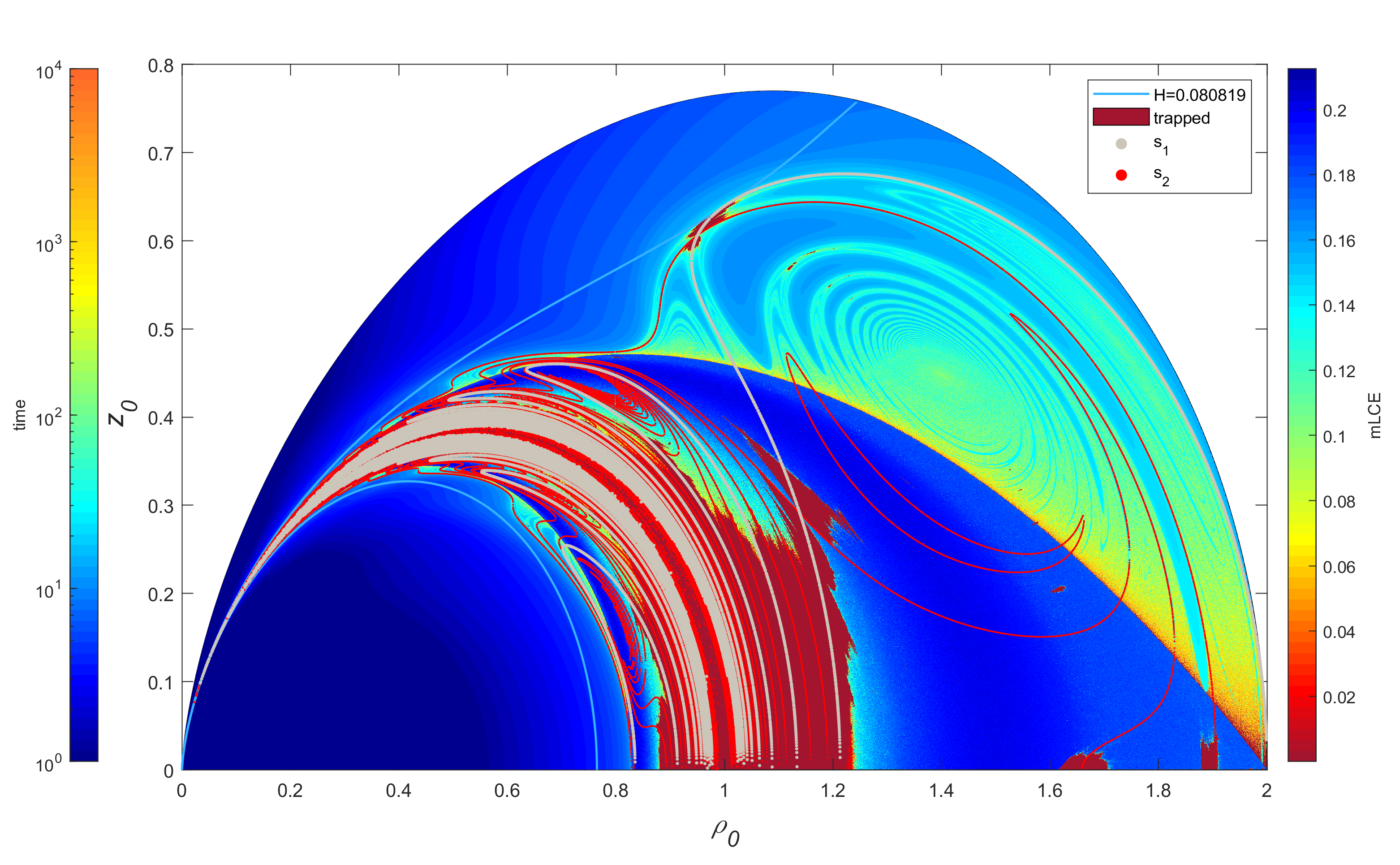}
\caption{\label{fig3} Same as Fig. \ref{fig1} with the open-path periodic orbits s2 and s1 indicated by red and grey lines, respectively.}
\end{figure*}

\begin{figure*}[htbp]
\includegraphics[width=1.8\columnwidth, angle=0.0]{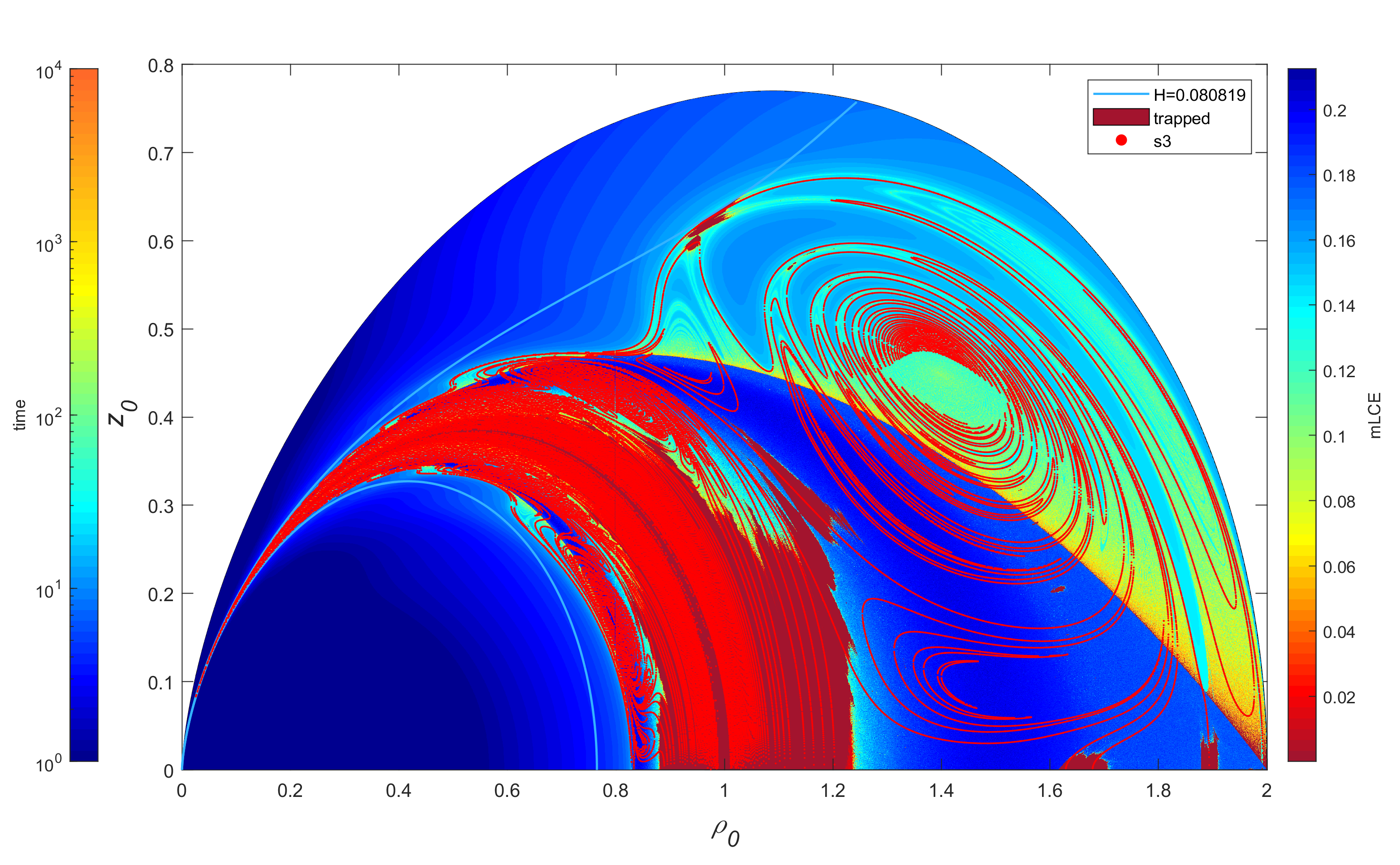}\\
\includegraphics[width=1.8\columnwidth, angle=0.0]{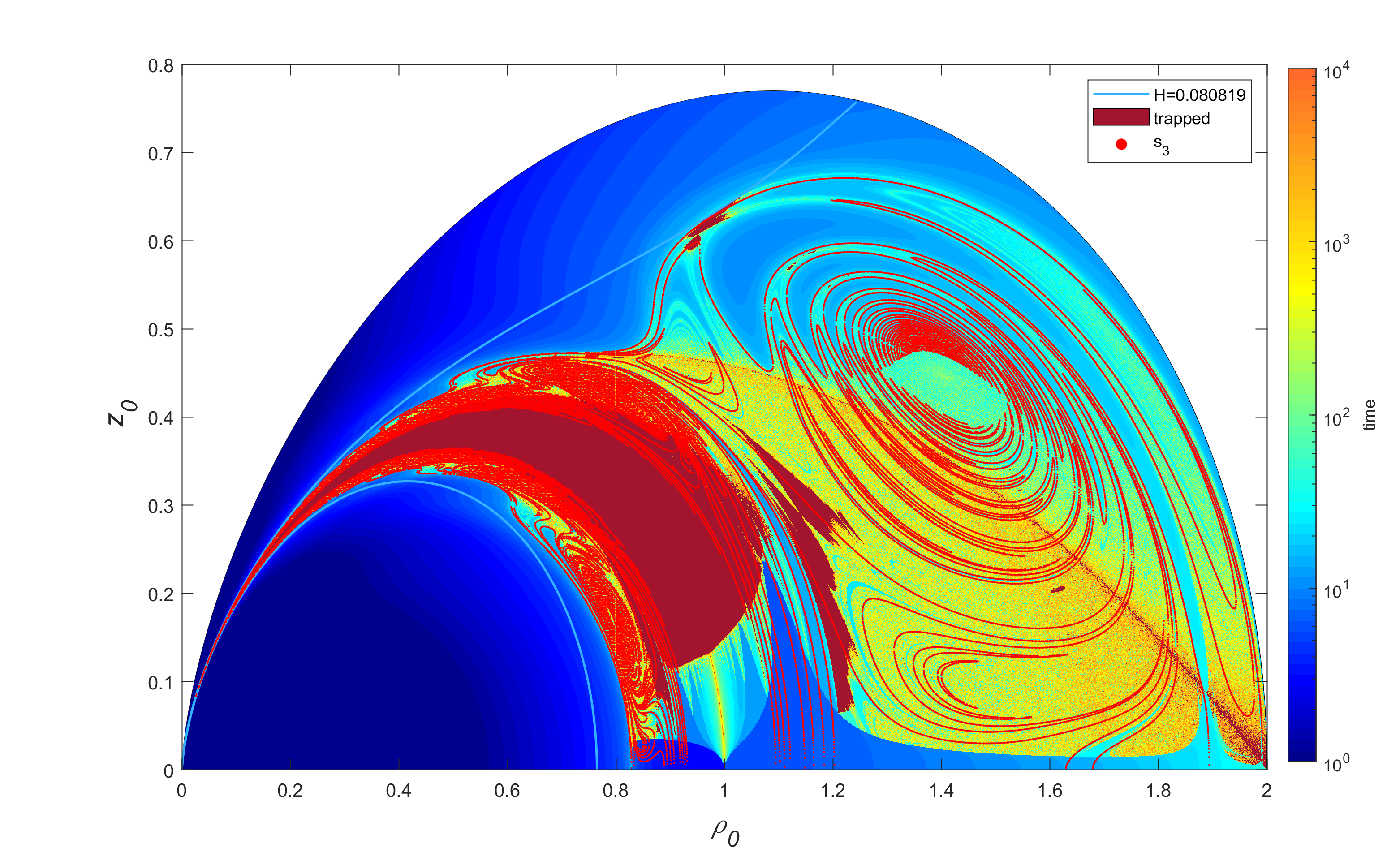}
\caption{\label{fig4} Upper: Same as Fig. \ref{fig1} with the open-path periodic orbits s3 indicated by red lines. Lower: Same as the upper panel except that the color for $H<1/32$ indicates the crossing time $t_{\rm cross}$ defined in the text.}
\end{figure*}

{\bf s1:} Given the general characteristics of symmetric open periodic orbits discussed in the previous section, one can readily identify classes of such orbits numerically. We first search for the principle families introduced by \citet{1978CeMec..17..215M}, members of which cross the equatorial plane once in each half a period. We search for such orbits by following the radial velocity $p_\rho$ at the first crossing of the equatorial plane of each orbit. Principle symmetric open-path periodic orbits have a null radial velocity at the first crossing of the equatorial plane. 

Figure \ref{fig2} shows families of such orbits (grey lines), which we labelled as "s1". They corresponds to families identified by \citet{1978CeMec..17..215M} with even indexes. However, in contrast to these families in the $[\rho,\ p_z]$ plane, each family locates on one side of the "thalweg", i.e., the bottom of the potential $V$. From high to low energies, we have $f_0$, $f_4$, ... and $f_2$, $f_6$, .. on the right and the left side of the "thalweg", respectively. The bulk of stable orbits associated with the low energy part of each family have similar features at the boundary, and $f_0$ goes through a prominent region of stable orbits with $H>1/32$, which agrees to the results of \citet{1978CeMec..17..233M}. Most stable orbits associated with the high-energy stable periodic orbits in the equatorial plane \citep{devogelaere1950} are not associated with families in s1, except for $f_2$, both of whose ends on the equatorial plane are stable.

{\bf s2:} With the principle families identified above, one can readily identify the second class of symmetric open period orbits, which cross the equatorial plane perpendicularly at the second crossing from the initial location. Numerically, when scanning the $(\rho_0, z_0)$ plane, we trace the radial velocity $p_\rho$ at the second crossing of the equatorial plane of each orbit. Periodic orbits in s2 have $p_\rho=0$ at the second crossing of the equatorial plane. The orbits obtained this way include those in s1, which have been identified above and can be readily separated from those in s2.

Figure \ref{fig3} shows s2 (red) together with s1 (grey). For each family in s1, there are 2 pairs of families in s2, with two families on each side of the thalweg, that go through the origin and terminate at 2 pairs of stable orbits in the equatorial plane at the other ends. Each pair of these stable orbits actually correspond to the two ends of one orbit in the equatorial plane. One of these families associated with $f_0$ in s2 crosses $f_0$ in the highest energy region of stable orbits. Another crosses $f_0$ at a lower energy. The orbit in s2 reduces to s1 at the crossing point of these two families. There are fine structures at the crossings between families in s1 and s2. The other two families on the other side of the thalweg do not go through families in s1.

The other families appear as loops and locate between these family pairs, and the number of loops increases quickly with the decrease of energy. These loops appear to trace spiral structures in the escape time and do not cross families in s1. Moreover, families in one class do not go through each other.

{\bf s3:} Similarly one can identify the third class of symmetric open periodic orbits that cross the equatorial plane perpendicularly at the third crossing from each of their initial location. Families in this class are shown in Figure \ref{fig4}. Similar to the class s2, for each family in s1, there are 4 pairs of families in s3 that go through the origin and terminate at 
4 pairs of stable periodic orbits in the equatorial plane at the other ends.
Two of these families go through $f_0$ in the highest energy region of stable orbits, while the other two cross $f_0$ at lower energies. Fractal structures are seen near these crossing points. The other four families at the inner side of the potential do not cross families in s1.

Compared to s2, there are many more loops between these families, that clearly trace spiral structures in the escape time and approach to their centers. To reveal similar spiral structures for $H<1/32$, one can introduce $\rho_{\rm max}=[1-(1-4(2H)^{1/2})^{1/2}]/[2(2H)^{1/2}]$ and $\rho_{\rm min}=[(1+4(2H)^{1/2})^{1/2}-1]/[2(2H)^{1/2}]$ that correspond to the maximum and the minimum distance to the origin for motion in the equatorial plane with an energy $H$, respectively. We define the crossing time $t_{\rm cross}$ as the first crossing of $\rho = \rho_{\rm max}-(\rho_{\rm max}-\rho_{\rm min})/50$ with a positive $p_\rho$. This crossing time is indicated in the lower panel of Figure \ref{fig4}. Similar to the escape time, initial conditions with the crossing time greater than $10^4$ are indicated with dark red color and are well correlated with stable orbits in the upper panel of Figure \ref{fig4}. There are spiral structures in the crossing time that are associated with family loops in s3. Similar to s2, the family loops do not appear to cross any other families, and families in s2 do not appear to cross families in s3 either, which is different from those in the [$\rho,\ p_z$] plane, where sub-branches of periodic orbits may cross each other \citep{1978CeMec..17..233M}.

\subsection{Classes of Asymmetric Open Periodic Orbits}

 \begin{figure*}[htbp]
\includegraphics[width=1.6\columnwidth, angle=0.0]{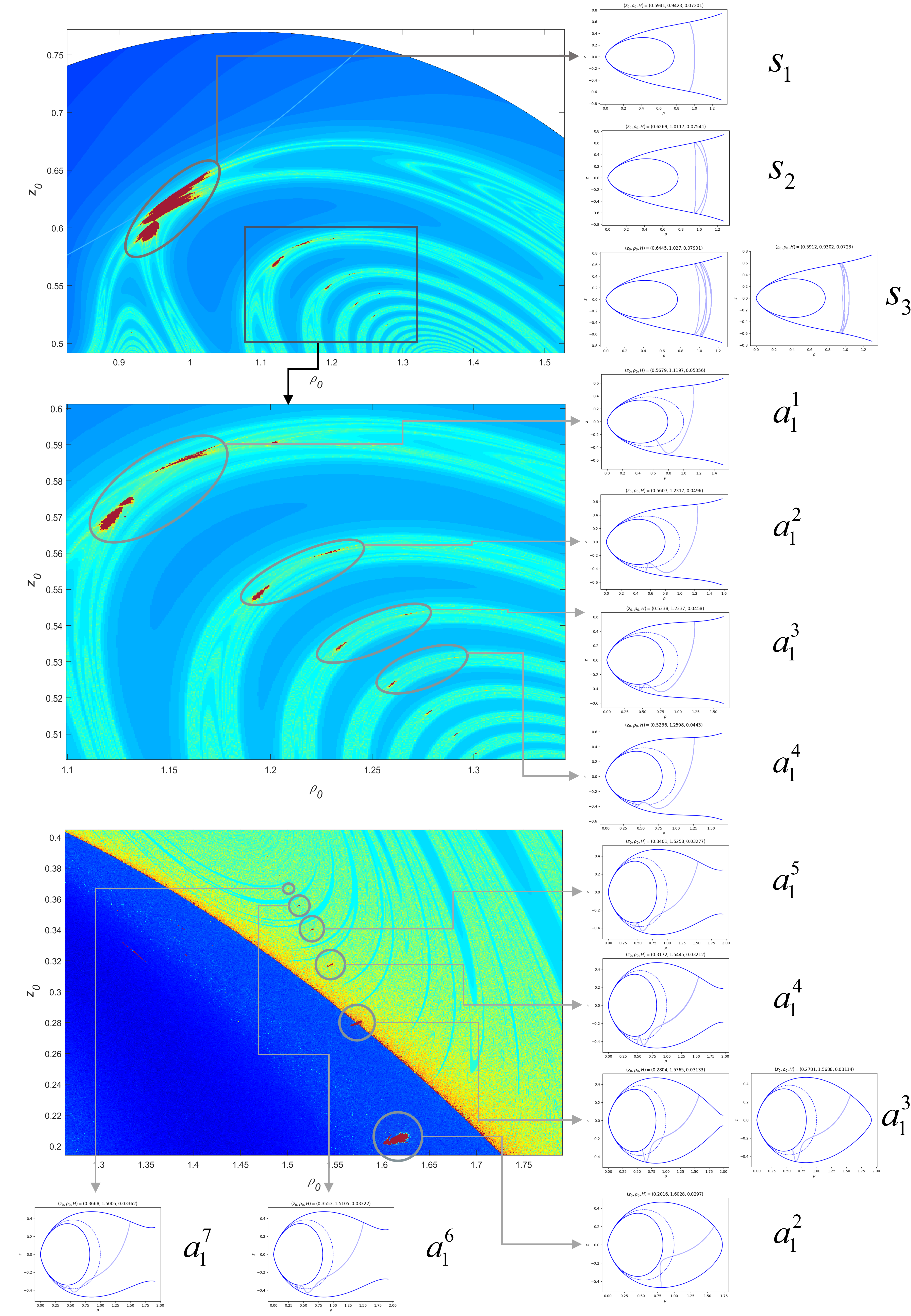}
\caption{\label{fig5} Examples of stable periodic orbits (thin solid lines) with the corresponding energy contours (thick solid lines) and the bottom of the potential $V$ (dotted) also shown.}
\end{figure*}

\begin{figure*}[htb]
\includegraphics[width=0.45\linewidth, angle=0.0]{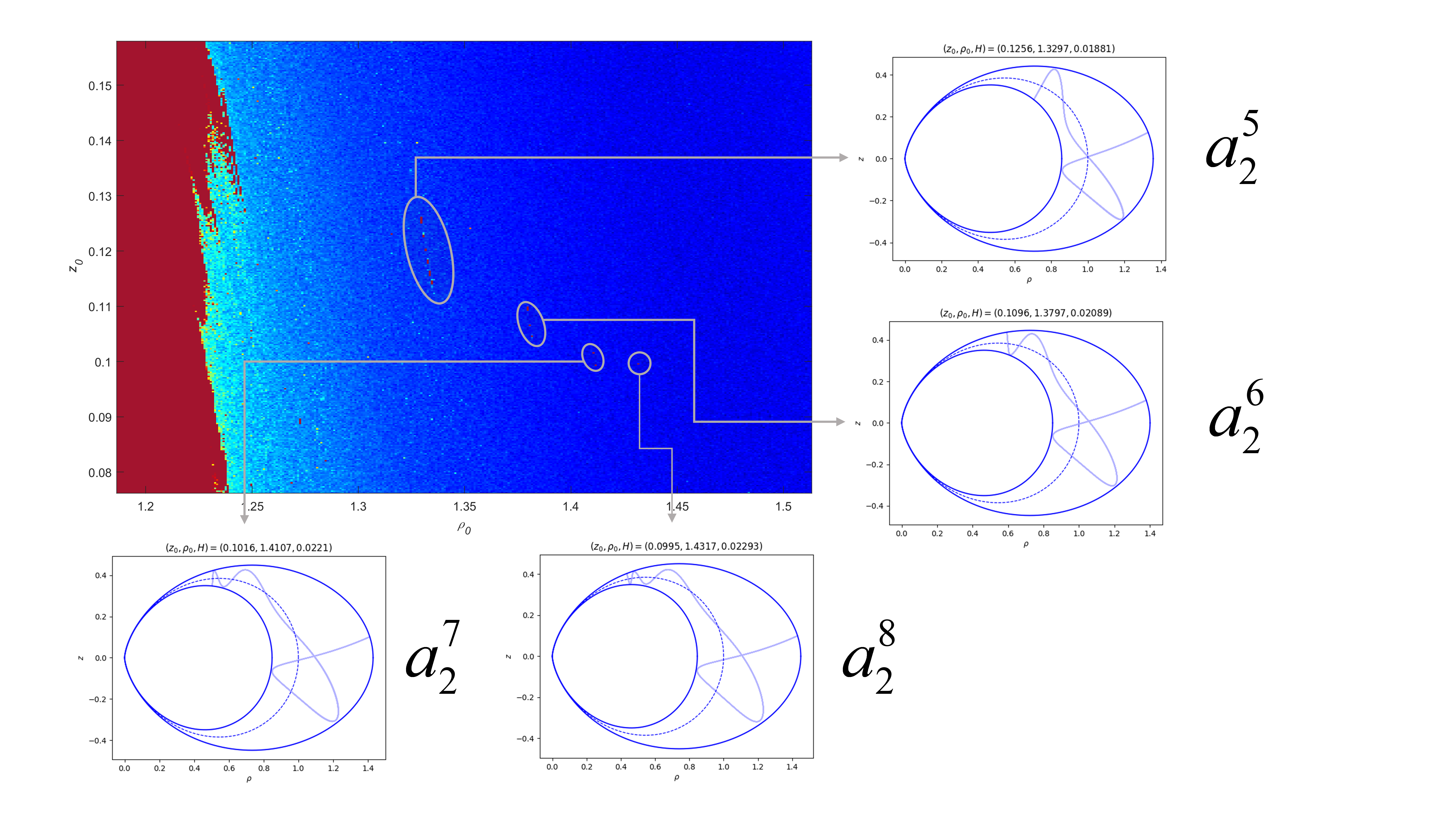}
\includegraphics[width=0.45\linewidth, height=0.27\linewidth, angle=0.0]{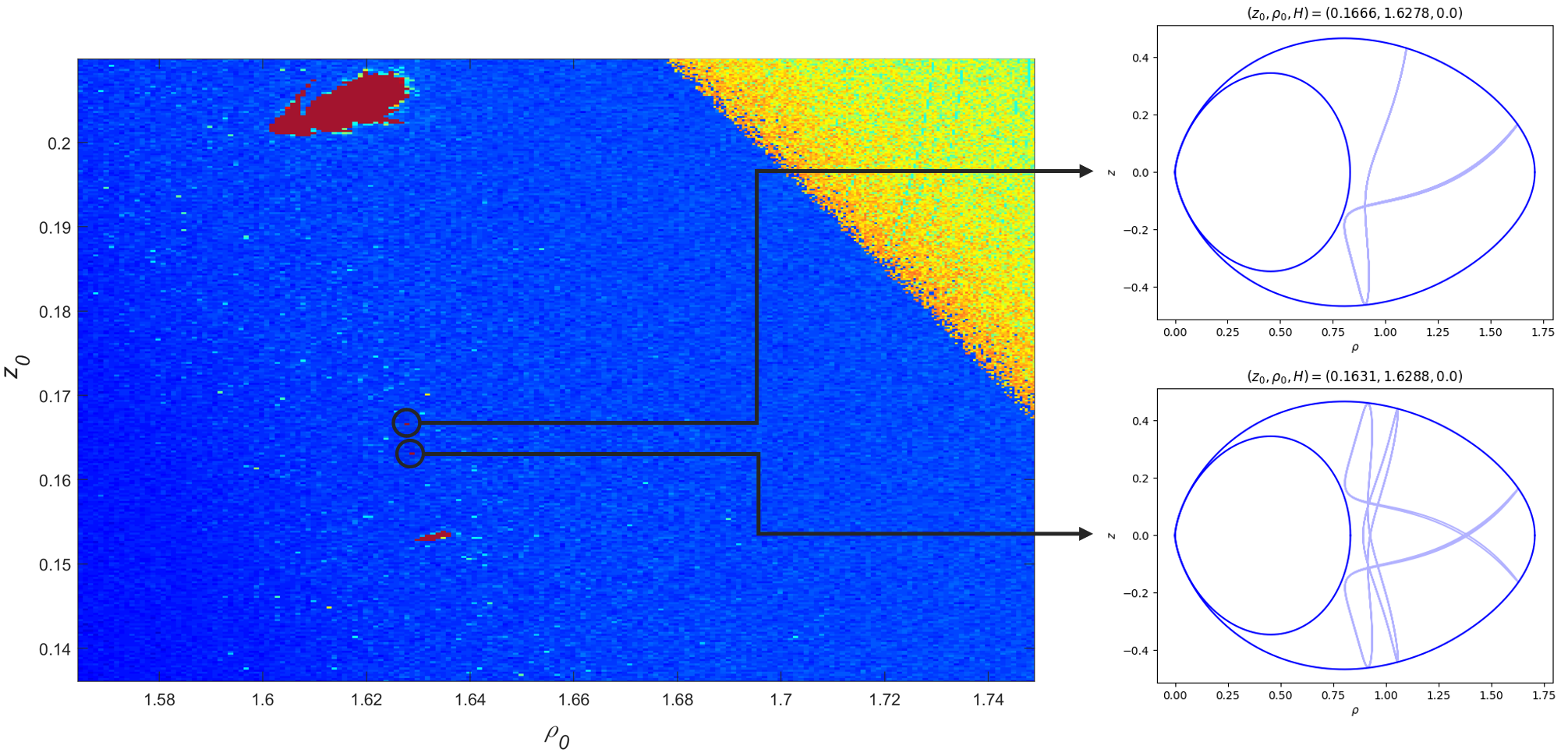}
\includegraphics[width=0.45\linewidth, angle=0.0]{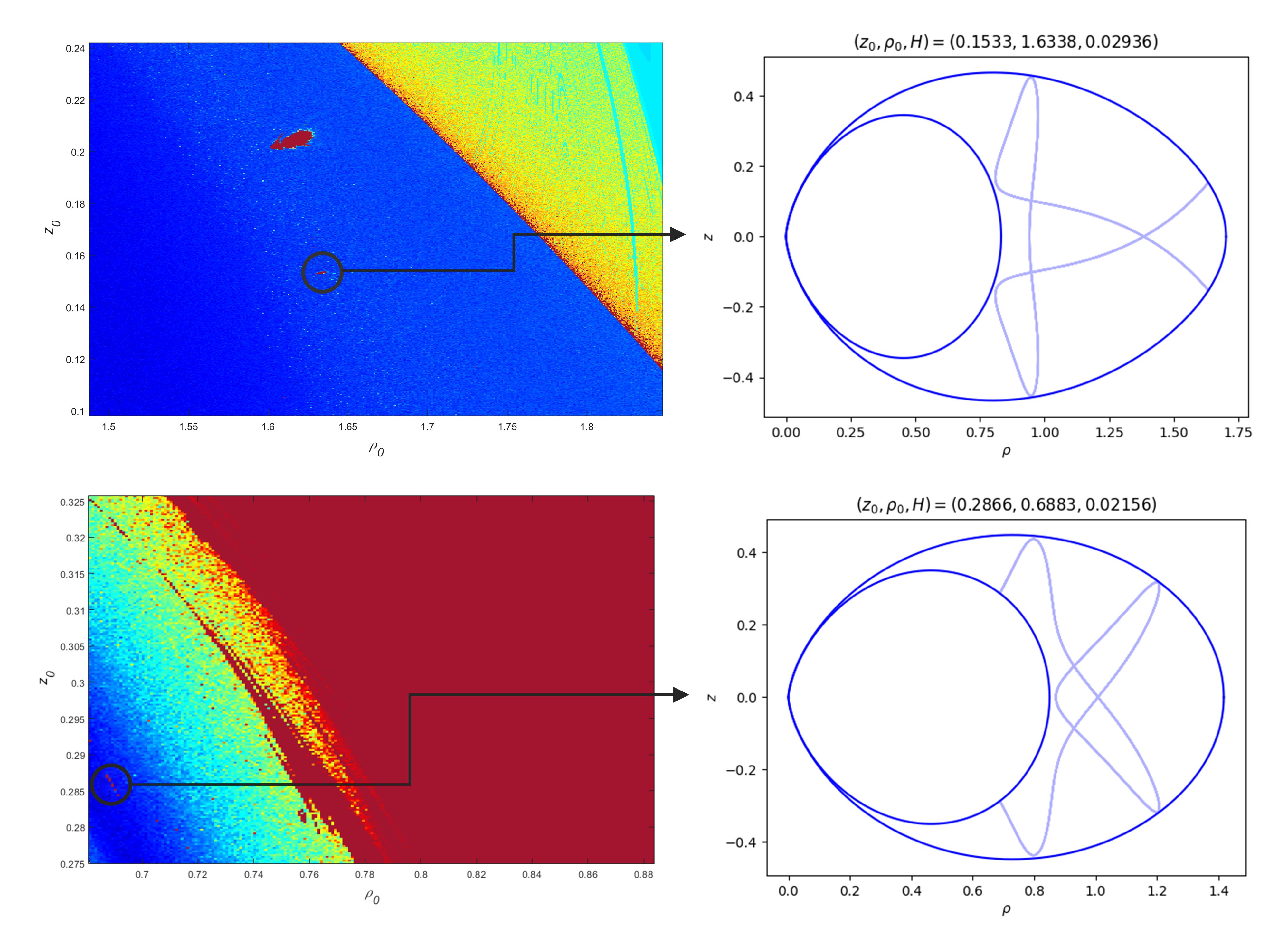}
\includegraphics[width=0.45\linewidth, height=0.34\linewidth, angle=0.0]{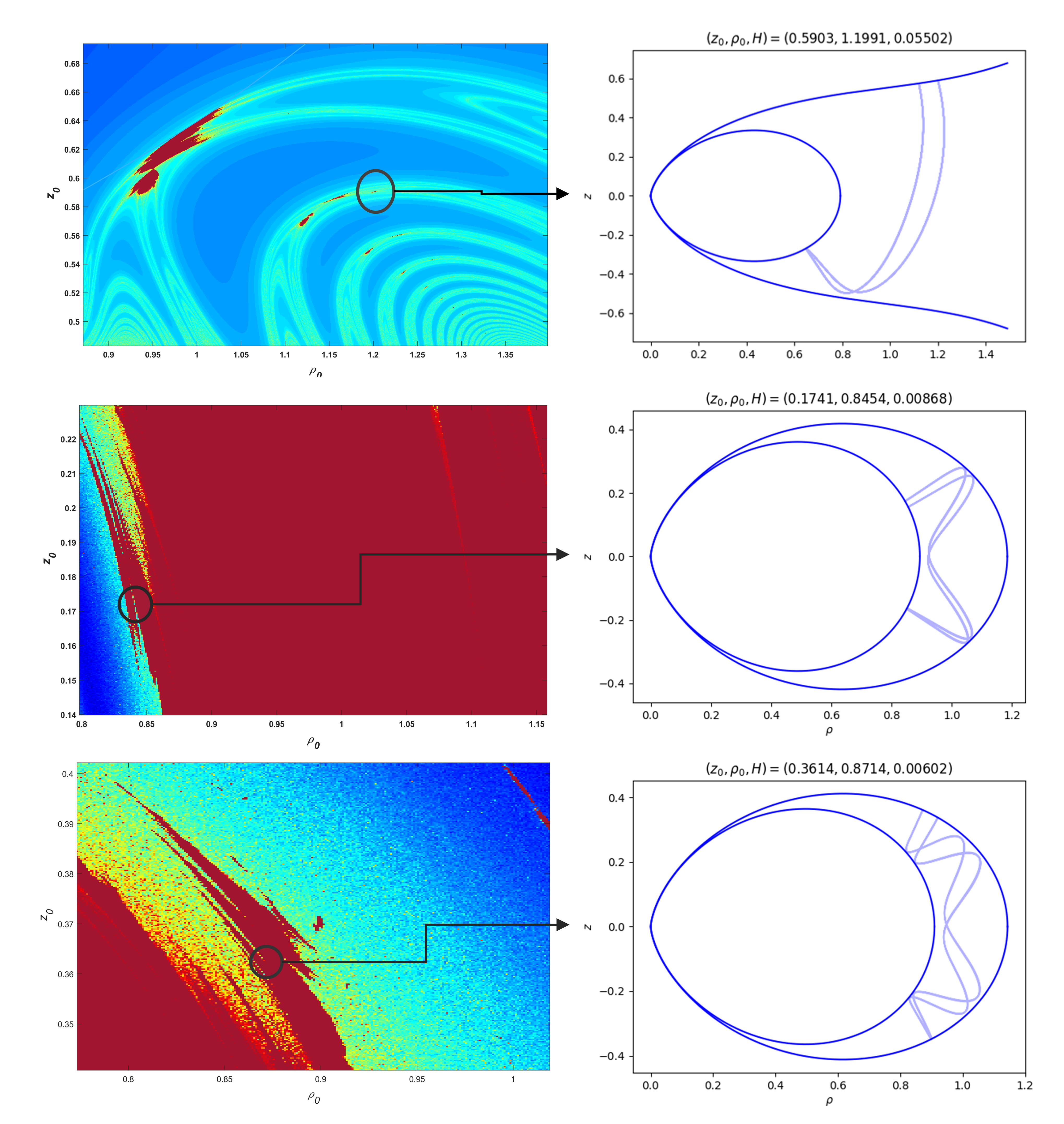}
\caption{\label{fig6} More examples of stable periodic orbits. }
\end{figure*}

Although the s3 class of symmetric open periodic orbits go through the most prominent regions of stable orbits, the chain of isolated regions of stable orbits at the turns of the main spiral structure are clearly not associated with s3 (See Fig. \ref{fig1}). We therefore randomly plot a periodic orbit from each of these regions to reveal their properties. Figure \ref{fig5} shows some of these orbits. The upper panel gives examples of symmetric periodic orbits in the highest energy region of stable orbits. Since all classes of symmetric periodic orbits go through this region, one can readily find orbits for each of these three classes presented above. Besides these trajectories, we also plot the corresponding energy contours to better demonstrate the orbital behavior.

The middle and lower panels show the upper and lower chains of islands of stable orbits. Based on the number of crossing through the equatorial plane and through the bottom of the potential $V$: $r=\cos^2\lambda$ in each half a period, we label these orbits with the letter "a" (for asymmetric) with the lower and upper indexes indicating these numbers, respectively. It is evident that orbits with an odd number of the upper index have a static point on each side of the thalweg while those with an even number of the upper index have their static points on one side of the thalweg. These two kinds of asymmetric orbits belong to different families. The two orbits $a_1^3$ in the lower panel show that for stable orbits near $H=1/32$, orbits with $H<1/32$ are very similar to those with $H>1/32$. The critical energy $H=1/32$ is not an important parameter determining properties of these orbits.  

On the other hand, along families of symmetric open periodic orbits, the number of crossing of the bottom of the potential may change \citep{1978CeMec..17..215M}. Moreover \citet{1977Ap&SS..48..471M} discovered a principle asymmetric class (surface) in the 3 dimensional sub-phase-space of $z=0$ that contains both open and close types of periodic orbits, and they both can oscillate around the bottom of the potential $V$ many times. They also found more than 20 regimes of stable orbits in this asymmetric class. In the 4 dimensional phase space, the trajectories of this asymmetric class occupy a 3 dimensional volume and its intersections with the plane of static points for the Meridian motion should be 1 dimensional. The stable asymmetric orbits in Figure \ref{fig5} all belong to the principle asymmetric class identified by \citet{1977Ap&SS..48..471M}!

Nevertheless, in the  plane of static points, they are divided into two families, one with their two ends on one side of the thalweg and the other with their two static ends on each side of the thalweg. Figure \ref{fig5} shows that these two families follow the main spiral structure and the number of oscillations around the bottom of the potential increases as these families approach the center of this spiral structure. Since these two families cross the equatorial plane once in each half a period, their two ends locate on two sides of the equatorial plane.

The principle asymmetric class terminates at the equatorial plane, the family with an odd upper index should terminate at the equatorial plane as well. The other ends of this family locate on the inner side of the thawhal (below the equatorial plane). As this family approaches the center of the spiral structure, its the other ends approach the origin.  This family therefore has another branch on the inner side of the potential. The family with an even upper indexes crosses $f_0$ in the low energy end and extends towards the origin. Since both ends of this family are on the outer side of the thalweg, there is a one to one map for points on both sides of the crossing point with $f_0$. As this family approaches the center of the main spiral structure, it approaches the origin at the other end (below the equatorial plane).

  \begin{figure*}[htbp]
\includegraphics[width=0.4\linewidth, angle=0.0]{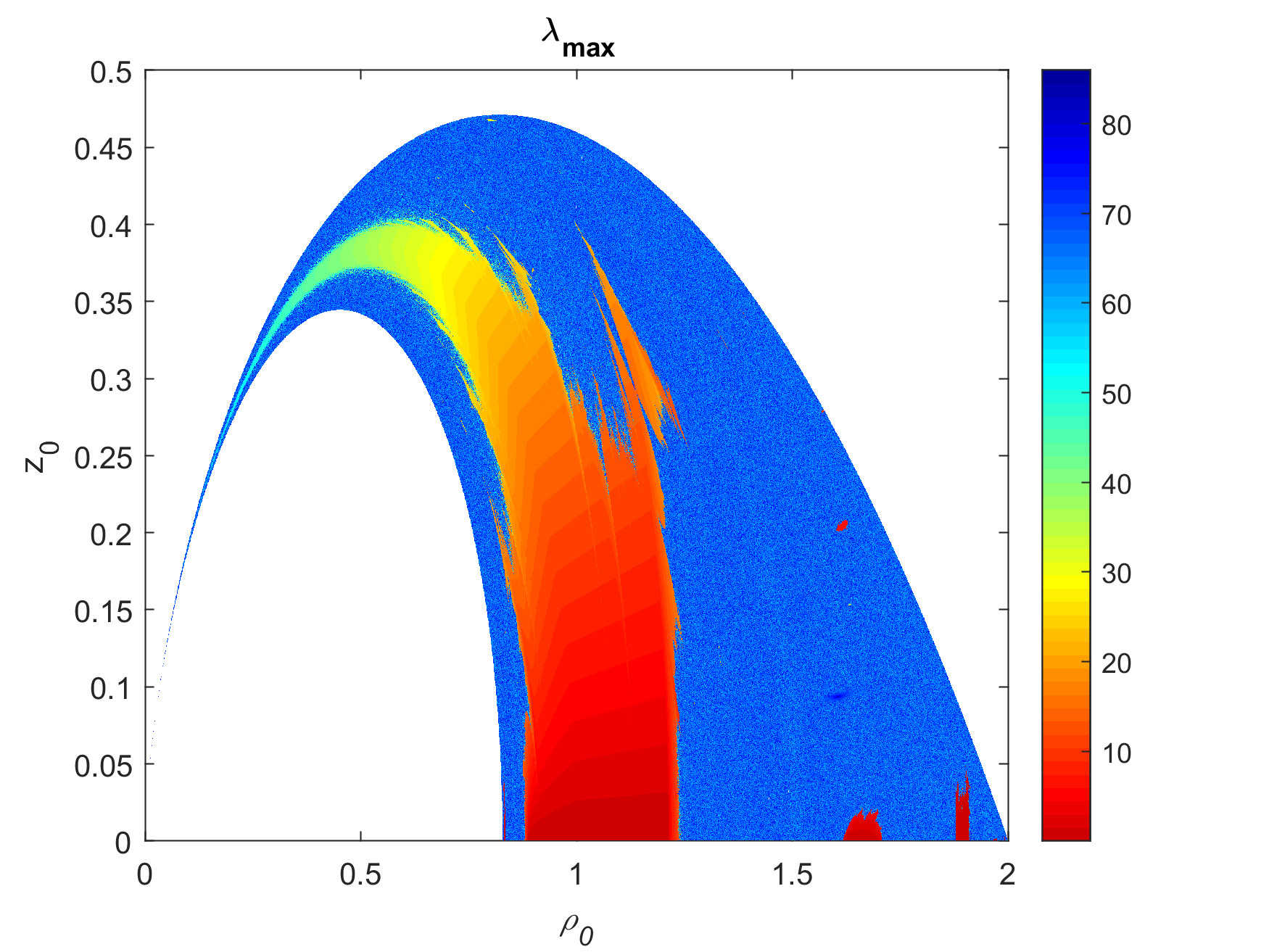}
\includegraphics[width=0.4\linewidth, angle=0.0]{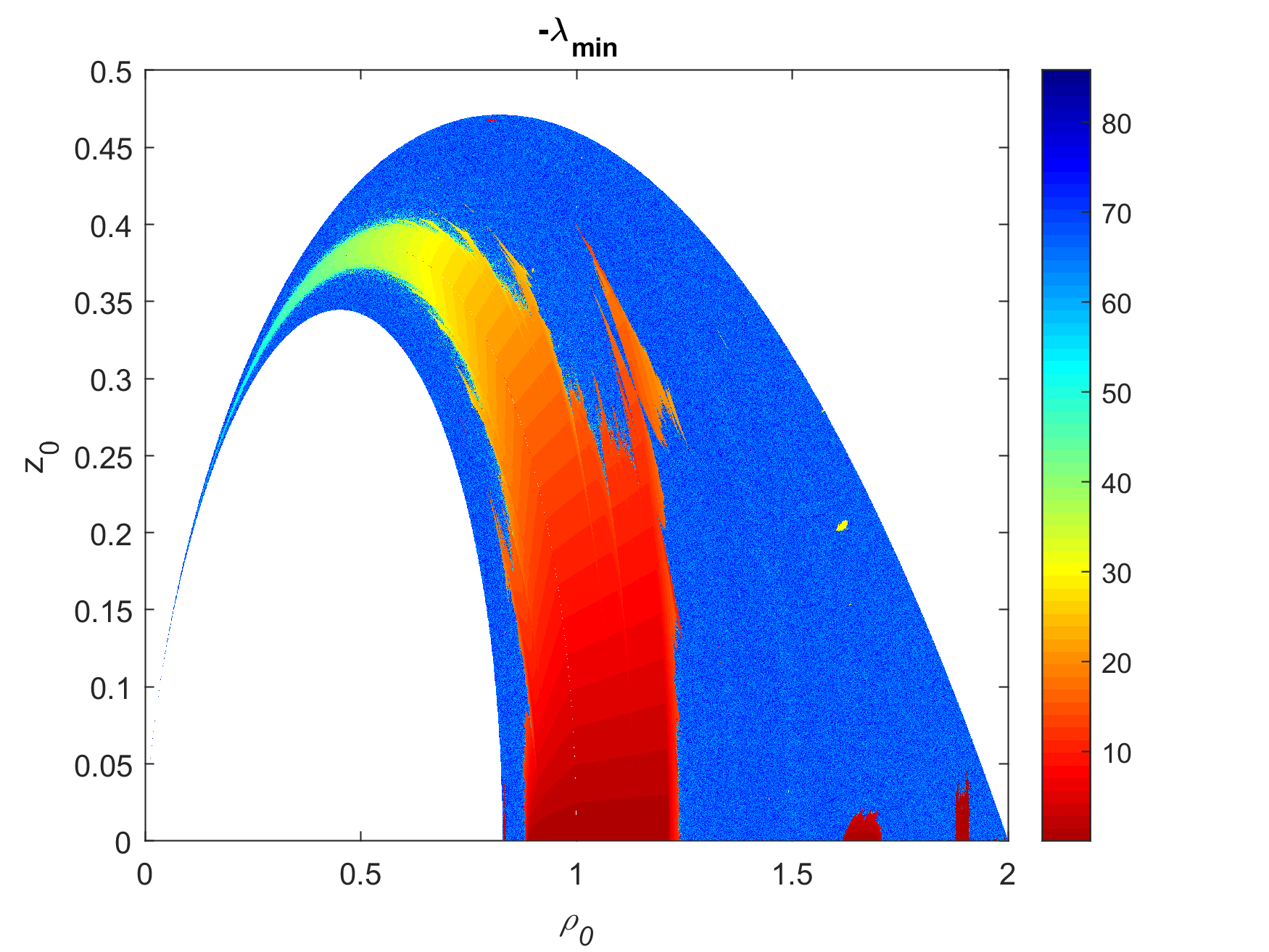}
\includegraphics[width=0.5\linewidth, angle=0.0]{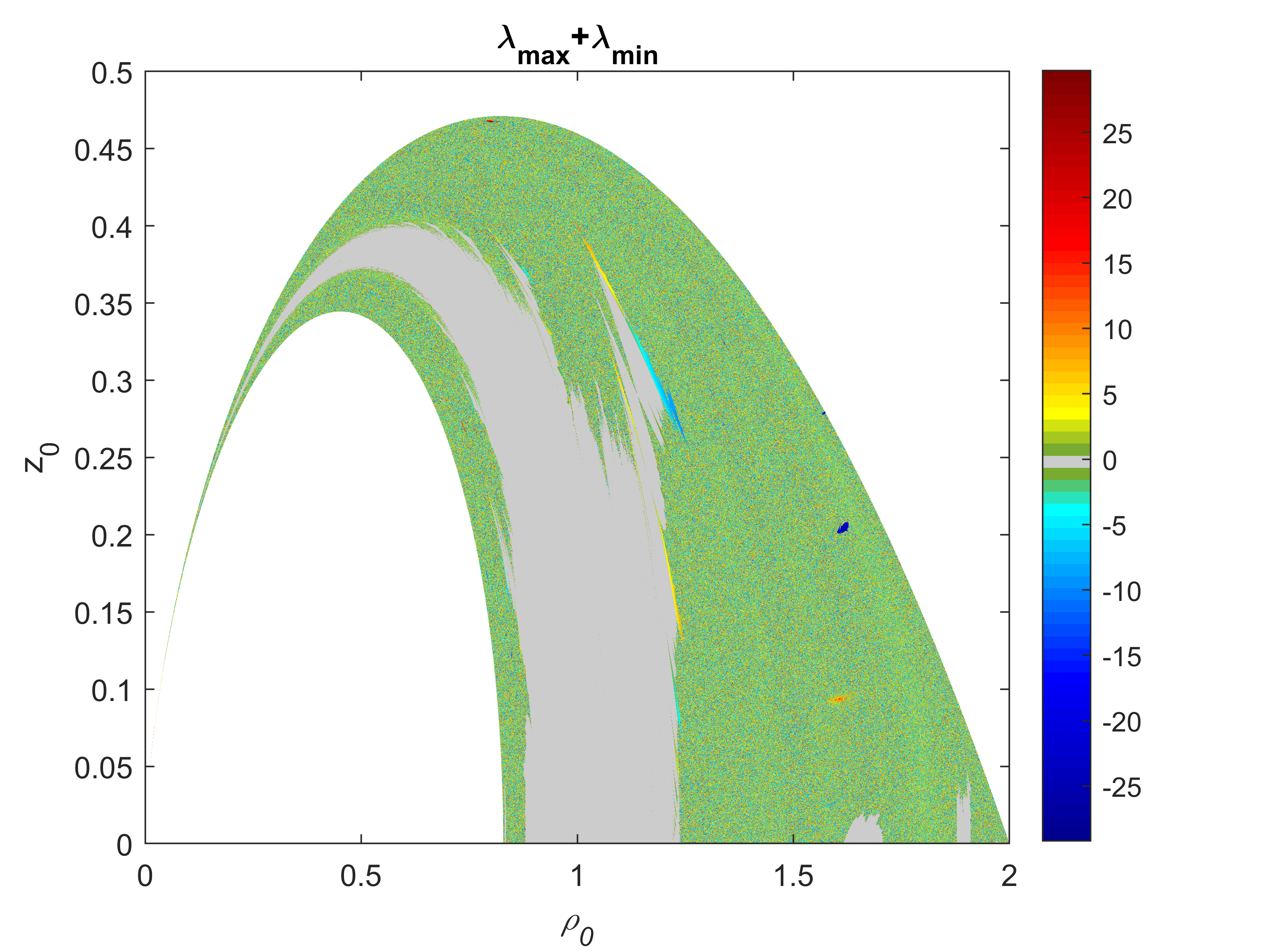}
\caption{\label{fig7} Upper: The maximum (left: $\lambda_{\rm max}$) and minimum (right: $\lambda_{\rm min}$) angles to the equatorial plane for $H<1/32$. Lower: the sum of the maximum and the minimum of $\lambda$.}
\end{figure*}

Besides the main spiral structure, there are many other smaller spiral structures. The upper left panel of Figure \ref{fig6} shows some stable asymmetric periodic orbits around a spiral structure below the main spiral structure and closer to the equatorial plane (See Fig. \ref{fig4}). Compared with the principle asymmetric class discussed above, it crosses the equatorial plane twice in half-a-period. Otherwise their behaviors are very similar to those around the main spiral structure. There therefore should be two families of asymmetric orbits around each spiral structure, and one of them crosses $f_0$ approaching the origin and the other terminates at the equatorial plane. Since both of their open ends are on one side of the equatorial plane, as the family with an even upper index crosses $f_0$, its two ends merge into one and another static point appears near one quarter of its period (i.e. half a period of the corresponding symmetric orbit). In contrast to the principle asymmetric class that bifurcates from $f_0$ with the same period, the period doubles as this class bifurcates from $f_0$.

\subsection{Other Regions of Stable Orbits}

We also notice that there are other tiny islands of stable orbits. In the upper right panel of Figure \ref{fig6}, the upper orbit is similar to the families in the upper left panel but belongs to the main spiral structure since it locates above a family in s3 that surrounds the main spiral structure. The lower orbit belongs to s3. It likely belongs to a branch bifurcating from the family of the principle asymmetric class that approaches the origin. The period also doubles at the point of bifurcation.  

The lower left panel of Figure \ref{fig6} shows two stable orbits that belong to families in s2 that approach the origin. The lower one actually belongs to a family that is associated with $f_2$. The top one in the lower right panel of Figure \ref{fig6} shows a stable asymmetric open orbit that appears to be bifurcating from the family in the principle asymmetric class that terminates at the equatorial plane. However, the top, middle, and bottom panel are associated with $f_0$, $f_2$, and $f_4$, respectively, confirming self-similarity among different families in s1.

\subsection{Asymmetry of Open Orbits}

Calculation of the mLCE is time consuming. There is a more efficient way to obtain region of stable asymmetric orbits. Given the asymmetry of the corresponding periodic orbits, one may just obtain the maximum and minimum values of $z/r$ for a given long enough computation time. Figure \ref{fig7} shows the maximum and minimum of $\lambda = \arcsin{z/r}$ and the sum of these two for a dimensionless computation time of $10^4$. The most prominent isolated region of stable asymmetric orbits near $(z=0.205\, \rho=1.62)$ can be readily seen. Moreover, we find some asymmetric stable orbits near the boundary of the region for the bulk of low energy stable orbits. It is interesting to note that the center of the spiral structure below the main spiral and closer to the equator is also highly asymmetric. 

\section{Conclusions and Discussions}
\label{conc}

In this paper, we systematically study orbits of charged particles in a dipole magnetic field with an azimuthal initial velocity that correspond to the open-path type orbits in the Meridian plane starting with a zero velocity. To understand the distribution of stable orbits, we identified three classes of open periodic orbits symmetric with respect to the equatorial plane. We found that orbits in the class s1 that cross the equatorial plane once in half a period locate on several separated families and there are compelling evidence of similarities among these families. We therefore focus on exploring stable orbits associated with the highest energy one $f_0$. 

We found that families in s2 and s3 either approach the origin with the other end terminating at the equatorial plane or have a loop shape. The former appears in pairs and one branch of these pairs cross $f_0$ at a stable orbit. These families separate the phase space into several regions with spiral structures in the escape and crossing times. The loops trace these spiral structures. Families in s2 and s3 do not seem to cross each other.

Although we did not identify classes of asymmetric open periodic orbits explicitly, by plotting orbits in isolated regions of stable orbits, we found that except for the highest energy one and those associated with stable orbits in the equatorial plane, the most prominent isolated regions of stable orbits are associated with a class of asymmetric open periodic orbits, that is the principle asymmetric family identified by \cite{1977Ap&SS..48..471M}. In the plane with zero meridian velocities, this class can be divided into two families: one crossing $f_0$ and approaching the origin with the two open ends of each orbit locating on either side of $f_0$ and the other one terminating at the equatorial plane with the other open end of each orbit in a branch on the other side of the bottom of the potential $V$. The other ends of these two families trace the main spiral structure and approach to its center. The numbers of oscillation around the bottom of the potential $V$ increase as one moves along these families toward the center of the main spiral structure with the other open end of the corresponding orbit approaching the origin. 

Loops in the symmetric open periodic orbit classes also approach the center of the main spiral structure with the increase in number of crossing of the equatorial plane in half a period. The center of the main spiral structure appears to singular in the sense that both symmetric periodic orbits oscillating around the equatorial plane and asymmetric periodic orbits oscillating around the bottom of the potential with their two open ends locating in regions near the origin and this central point approach to this same point with the increase in the number of oscillations. The nature of the orbit at the center of spiral structure is ambiguous. One can imagine that similar families of asymmetric orbits exist around each spiral structures in the escape and crossing times. There are therefore many singular centers. Further explorations are warranted.

We found evidence that, besides the principle asymmetric family bifurcating from $f_0$ \citep{1978CeMec..17..233M} and tracing the main spiral structure, there are many other asymmetric families bifurcating from $f_0$ and tracing other spiral structures.  There are also branches bifurcating from the asymmetric families tracing these spiral structures. Such kinds of bifurcation exist in all families in s1.

\begin{acknowledgments}
We appreciate helpful discussions with Ms. Yilin Yang. This work is partially supported by the National Key R\&D program of China under the grant Nos. G2021166002L, and 2018YFA0404203, NSFC grants U1931204, 12147208, 11947404 and 11761131007, DFG Sino-German Collaboration Project Nos. BU 777/15-1 and MU 4255/1-1, Department of Science and Technology of Sichuan Province No.2020YFSY0016, and by the SRTP program of the school of physical science and technology, Southwest Jiaotong University No. 202110613088. 
\end{acknowledgments}

\section*{Data availability}

The data that support the findings of this study are available from the corresponding author upon reasonable request.



\bibliography{sorsamp} 

\end{document}